\documentclass[sii]{ipart}

\RequirePackage[OT1]{fontenc}
\RequirePackage{amsthm, amsmath, amssymb}
\RequirePackage[numbers,square]{natbib}
\RequirePackage[colorlinks]{hyperref}
\usepackage{bm}
\usepackage{mathtools}
\usepackage{mathrsfs}
\usepackage{threeparttable}

\pubyear{0000}
\volume{0}
\issue{0}
\firstpage{1}
\lastpage{14}

\startlocaldefs
\theoremstyle{plain}

\makeatletter

\newcommand{\Rmnum}[1]{\expandafter\@slowromancap\romannumeral #1@}
\makeatother

\newtheorem{theorem}{Theorem}

\newcommand{\bZ}{{\bf Z}}
\newcommand{\bz}{{\bf z}}

\newcommand{\bu}{{\bf u}}
\newcommand{\bbeta}{\boldsymbol{\beta}}

\newcommand{\T}{\!\top\!}

\def\E{ {\sf E}}

\endlocaldefs

\allowdisplaybreaks

\begin{document}

\begin{frontmatter}
\title{Proportional Mean Residual Life Model with Censored Survival Data under Case-cohort Design}
\runtitle{Proportional Mean Residual Life Model with Censored Survival Data under Case-cohort Design}

\begin{aug}
\author{\fnms{Huijuan}
\snm{Ma}\ead[label=e1]{mahuijuan2015@163.com}},
\address{Department of Biostatistics and Bioinformatics,\\
Emory University,\\
Atlanta, GA, 30322, USA\\
\printead{e1}}

\author{\fnms{Jianhua} \snm{Shi} \thanksref{t1} \ead[label=e2]{v0085@126.com}}
\address{School of Mathematics and Statistics,\\
Minnan Normal University,\\
Zhangzhou, Fujian, 363000, China\\
\printead{e2}}
\and
\author{\fnms{Yong} \snm{Zhou}
\ead[label=e3]{yzhou@amss.ac.cn}}
\address{Institute of Statistics and Interdisciplinary Sciences \\
and School of Statistics,\\
Faculty of Economics and Management, \\
East China Normal University,\\
Shanghai, 200241, China\\
\printead{e3}}\\
\thankstext{t1}{The corresponding author.}
\runauthor{Ma et al.}
\end{aug}

\received{\sday{30} \smonth{7} \syear{2017}}

\begin{abstract}
Proportional mean residual life model is studied for analysing survival data
from the case-cohort design.
To simultaneously estimate the
regression parameters and the baseline mean residual life function, weighted estimating equations based on an inverse selection probability are proposed.
The resulting regression coefficients estimates are shown to be consistent and
asymptotically normal with easily estimated variance-covariance.
Simulation studies show that the proposed estimators perform very well.
An application to a real dataset from the South Welsh nickel refiners study
is also given to illustrate the methodology.
\end{abstract}


\end{frontmatter}

\section{Introduction}

When studying the natural history of a event, such as the fields of survival analysis, medical study, actuarial science
and reliability research, the residual lifetime is often regarded as a crucial index for investigators to make decisions.
The mean residual life function (MRLF) is one of the most important quantitative measures for the residual lifetimes that can
describe the characteristics of the residual life time more directly.
The MRLF for a nonnegative survival time $T$ with finite expectation at time $t \geq 0$
is defined as $m(t)=\E (T-t|T>t)$.
It is often of interest to analyse the mean residual life function in many applications.
For example, a driver may be interested in knowing how much longer his or her
car can be used, given that the car has worked normally for $t$ years.
Many early literatures on MRLF studied its probability behaviours, statistical inference on testing procedures
and the estimation in homogeneous cases.
Apparently, the MRLF may vary due to different covariates.
To quantify and summarize the association between the MRLF and its associated covariates,
extensive regression models are explored.
\cite{Oakes1990} originally proposed the proportional mean residual life model, which has been studied by many authors later.
The proportional mean residual life model, or the Oakes--Dasu model, is specified by
\begin{eqnarray}
m( t | \bZ ) = m_0(t) \exp(\bbeta^{\T} \bZ),
\end{eqnarray}
where $m( t | \bZ )$ is the mean residual life with the $p$-vector covariate $\bZ$,
$\bbeta$ is the usual regression $p$-parameter vector, and $m_0(t)$ is an unknown baseline mean residual life.
\cite{Maguluri1994} developed estimation procedures for regression coefficients mainly for
uncensored survival data, which was later modified to accommodate right censoring setting in \cite{Chen2005b}.
\cite{Chen2005} used counting process theory to develop another semiparametric
inference procedures for the proportional mean residual life model.
\cite{Chen2006} and \cite{Chen2007} proposed the additive mean residual life model
and discussed various estimation methodologies with or without right censoring.
\cite{sun2009class} proposed a more general family of transformed mean residual life model, including
the proportional mean residual life model and the additive mean residual life model as special cases.

However, the above methods for mean residual life models are not suitable when some covariates are missing.
In large cohort studies, the major effort and cost arise from the assembling and analysing of covariate measurement.
When the disease rate is low, assembling all covariates for every subject may become redundant and expensive.
\cite{prentice1986case} proposed case-cohort design
to provide a cost effective way of
conducting such cohort studies.
Under this design, a random sample from the entire cohort is selected, named the subcohort.
Covariate information is collected only for the subjects in the subcohort and all the cases
who experience the event of interest.
After the landmark article of \cite{prentice1986case}, the case-cohort
design has been extensively studied in many statistical literatures.
Standard analysis of the case-cohort design are conducted using the Cox proportional hazards model \citep{Cox1972}.
 For example, a pseudo-likelihood procedure proposed by \cite{prentice1986case} was later elaborated by \cite{Self1988}, \cite{Lin1993} and \cite{Ma2017}.
Several authors studied other regression
models such as the additive hazards model \citep{Kulich2000, Ma2015, Li2017}, the proportional odds model
\citep{Chen2001a} and the semiparametric transformation regression model \citep{Chen2001b, Kong2004, Lu2006, Chen2009, Ma2016}.
\cite{Borgan2000}, \cite{Kulich2004} and \cite{Breslow2007}
among others, extended the classical case-cohort design to more complex sampling schemes.
Besides, \cite{Zheng2013} and \cite{Fan2018} conducted quantile regression analysis of case-cohort data.
All these models may be adopted to indirectly make statistical inference for the mean residual lifetime.
But they are relatively cumbersome and not straightforward to measure the residual life.
Further, the mean residual life function is appealing to understand for practical use, and it provides an alternative to the hazard function.
Consequently, improving statistical methods for mean residual life models are needed under the case-cohort design.

Just as mentioned before, existing methods for mean residual life models are
used for cohort data with complete covariate information, they are not suitable for the case-cohort data.
To the best of our knowledge, there have been no study about mean residual life model under case-cohort design.
In this paper, we focus on proportional mean residual life model
for the analysis of case-cohort data.
Our research is initially  motivated by a nickel refiners study in the South Welsh where the refiners are interested in knowing how long they can still survive given  their
current situation. Thus, the mean residual life model is an informative choice.
Further, the event rate for this study is quite low and hence the case-cohort design is preferred.
Our approach is motivated by \cite{Chen2005}, which made use of
the counting process theory in constructing some estimating equations
and does not require estimating or modelling the distribution of
censoring.
The main difficulty here is that some covariates are missing and the subjects whether they should be selected in the subcohort are not
independent with each other.

The remainder of this paper is organized as follows.
In Section 2, several new weighted estimating equations are proposed for simultaneous estimation of the
regression parameters and the baseline mean residual life function.
Some large sample properties of the resulting regression coefficients estimates are also given in this Section.
Section 3 is devoted to simulation studies to examine the finite sample properties of the regression parameter estimators.
In Section 4, a real dataset from the South Welsh nickel refiners study is used to illustrate the proposed estimating procedures.
Section 5 contains some concluding remarks and the outline of the proofs is provided in the Appendix finally.

\section{Estimating Equations and Theoretical Results}

The failure time and potential censoring time are denoted as $T$ and $C$, respectively, which are assumed to be independent given the $p\times 1$ covariate vector $\bZ$.
Let $\tilde{T} = \min(T, C)$ and $\delta = I(T \leq C)$, then the usual counting process and the at-risk process at time $t$ can be defined as $N(t)=\delta I(\tilde{T}\leq t)$ and $Y(t)=I(\tilde{T}\geq t)$, respectively.
Complete data on a sample of $n$ individuals are modeled as $n$ independent and identically distributed random vectors
$(\tilde{T}_i,\delta_i,\textbf{Z}_i)$, where $\tilde{T}_i=\min(T_i, C_i)$ and $\delta_i=I(T_i \leq C_i)$ for $i=1,2,\cdots,n$.
Consider the filtration defined by $\mathcal{F}_t = \sigma \{N_i(u), Y_i(u), \bZ_i: 0 \leq u \leq t, i=1,2,\cdots, n\}$,
then
$M_i(t; \bbeta_*, m_*) = N_i(t) - \int_0^t Y_i(s) d \Lambda_i(s;  \bbeta_*, m_* )$
are martingales with respect to $\mathcal{F}_t$, where
$\Lambda_i(\cdot)$ denotes the usual cumulative hazard function for subject $i$,
$\bbeta_*$ and $m_*(\cdot)$ are the true values of $\bbeta$ and $m_0(\cdot)$, respectively.
The martingale properties of $M_i(\cdot)$ implies $\E [  d M_i(t; \bbeta_*, m_*) ] = 0$ for $i=1,\ldots,n$.
Furthermore,
\begin{eqnarray}
& & \E \left[ m_*(t) d M_i(t; \bbeta_*, m_*) \right]  \nonumber \\
&=& \E \left[ m_*(t) d N_i(t) - m_*(t)Y_i(t) d \Lambda_i(t;  \bbeta_*, m_* )  \right]  \nonumber \\
&=& \E \left[ m_*(t) d N_i(t) - Y_i(t) \{ \exp( -\bbeta_*^{\T} \bZ ) d t + d m_*(t ) \} \right] \nonumber \\
&=& 0.
\end{eqnarray}
For simplicity, we assume
that $ 0 < \tau = \inf\{t: \Pr( \tilde{T} > t) = 0 \} < \infty.$
When the gathered data was complete, \cite{Chen2005} proposed the following two estimating equations to estimate $(\boldsymbol{\beta}_*, m_*(t))$:
\begin{eqnarray}
& & \sum_{i=1}^n \big[ m_0(t)dN_i(t)- Y_i(t) \big\{\exp(-\boldsymbol{\beta}^{\T}\textbf{Z}_i) dt  \nonumber  \\
& & ~~~~~~+ d m_0(t) \big\} \big] = 0  ~~~(0\leq t\leq \tau),  \label{eq11}   \\
& & \sum_{i=1}^n \int_0^\tau \bZ_i \big[ m_0(t)dN_i(t)- Y_i(t) \big\{ \exp(-\boldsymbol{\beta}^{\T} \bZ_i) dt  \nonumber  \\
& & ~~~~~~+ d m_0(t) \big\} \big] = 0. \label{eq12}
\end{eqnarray}

Under the case-cohort design, since $\bZ_i$ is not observed for all subjects, the estimating equations
(\ref{eq11}) and (\ref{eq12}) based on the entire cohort data are no longer available.
In this paper, we assume that the subcohort with fixed size $\tilde{n}$ is drawn from the entire cohort
by the simple random sampling.
Let $\xi_i$ be the subcohort indicator, taking the value $1$ or $0$, whether the subject is included in the subcohort or not.
Hence the data can be
summarised as $\{(\tilde{T}_i,\delta_i,\xi_i,[\delta_i+(1-\delta_i)\xi_i] \bZ_i),i=1,2,\cdots,n\}$, which means that
$(\tilde{T}_i, \delta_i)$ are available for all individuals in the entire cohort, and $\mathbf{Z}_i$ only for
subjects in the subcohort with $\xi_i=1$, and all the cases outside the subcohort with $\delta_i=1$ and $\xi_i=0$.  Here $\xi_i$ is independent of
$(\tilde{T}_i,\delta_i, \bZ_i), i=1,2,\cdots,n$, while the $\xi_i's$ are dependent because of the sampling without replacement.
Similar to \cite{Kong2004} and \cite{Lu2006},
for each individual in the full cohort,
we define a weight $\pi_i=\delta_i+(1-\delta_i)\xi_i/\hat{p}$ by the idea of the inverse selection
probabilities, where $\hat{p}=\tilde{n}/n$.

Now we propose the estimator by two steps in the following.

First, we develop a new estimator for $m_0(t)$ as if the true regression coefficients $\bbeta$ have been known.
We propose the estimating equation by incorporating the weight $\pi_i$,
\begin{eqnarray}
 & & \sum_{i=1}^n \pi_i \big[ m_0(t)dN_i(t)-  Y_i(t) \big\{\exp(-\boldsymbol{\beta}^{\T} \bZ_i) dt \nonumber \\
 & & ~~~~~ + d m_0(t) \big\} \big]=0  \label{eq113}  \\ [1mm]
 & \Longleftrightarrow &  \left\{ \frac{\sum_{i=1}^n   \pi_i d N_i(t)}{\sum_{i=1}^n \pi_i Y_i(t)} \right\} m_0(t)- d m_0(t) \nonumber \\
 & & = \frac{\sum_{i=1}^n \pi_i Y_i(t)\exp(-\boldsymbol{\beta}^{\T} \bZ_i)}{\sum_{i=1}^n \pi_i Y_i(t)} d t.
\label{eq13}
\end{eqnarray}
Equation (\ref{eq13}) is in fact
a first-order linear ordinary differential equation about $m_0(t)$, which possesses a closed form solution
\begin{eqnarray}
\hat{m}_0(t) &\doteq& \hat{m}_0(t;\boldsymbol{\beta})  \nonumber \\
&=& \frac{1}{ S_n(t) } \int_t^\tau S_n(u) B_n(u;\boldsymbol{\beta} )d u,      \label{mceq7}
\end{eqnarray}
where
\begin{eqnarray*}
S_n(t) &=& \exp \left\{ -\int_0^t \frac{\sum_{i=1}^n  \pi_i d N_i(u)}{\sum_{i=1}^n \pi_i Y_i(u)} \right\} ~~\textrm{and} ~~ \\
B_n(t; \bbeta) &=& \frac{\sum_{i=1}^n \pi_i Y_i(t)\exp(-\bbeta^{\T}\bZ_i)}{\sum_{i=1}^n \pi_i Y_i(t)} .
\end{eqnarray*}
Based on the mean-zero process $ N_i(t)-\int_0^t \pi_i Y_i(u) d \Lambda_i(u; \bbeta_*, m_*)$, it is easy to show that $S_n(t)$
is a consistent estimator of the survival function for the failure time $T$.
We can see that (\ref{mceq7}) owns similar formula as $m(t) = \int_t^{\tau} S(u) d u /S(t)$ with $S_n(t)$ as an unbiased estimator of $S(t)$, while the additional term $B_n(u;\boldsymbol{\beta} )$ that involves $\bbeta$ is caused by the proportional assumption for the mean residual life model.

Next, we propose the estimating equations to estimate $\bbeta_*$.
Define
\begin{eqnarray} \label{eq002}
& & U\{ \bbeta, m_0(\cdot) \} =\frac{1}{n} \sum_{i=1}^n  \int_0^\tau
\pi_i \bZ_i \big[ m_0(t) dN_i(t)  \nonumber \\
& &~~~~~ - Y_i(t) \big\{\exp(-\bbeta^{\T} \bZ_i) dt + d m_0(t) \big\} \big].
\end{eqnarray}
Note that
\begin{eqnarray}
 U\{ \bbeta_*, m_*(\cdot) \}
&=&  \frac{1}{n} \sum_{i=1}^n  \int_0^\tau \pi_i \bZ_i \big[ m_*(t) dN_i(t) \nonumber \\
& & - Y_i(t) \big\{\exp(-\bbeta_*^{\T} \bZ_i) dt + d m_*(t) \big\} \big]  \label{eqap3} \\
&=&  \frac{1}{n} \sum_{i=1}^n \int_0^\tau \pi_i \bZ_i  m_*(t) d M_i(t; \bbeta_*, m_*).  \label{eqap4}
\end{eqnarray}
are mean-zero.
 To obtain a consistent estimator for $\boldsymbol{\beta}_*,$ we replace $m_0(t)$ with
 $\hat{m}_0(t;\boldsymbol{\beta})$ in $U\{ \bbeta, m_0(\cdot) \}$, and define $ \bar{\textbf{Z}}(t)= \sum_{i=1}^n \pi_i \textbf{Z}_i Y_i(t)/\sum_{i=1}^n \pi_i Y_i(t)$,
 then the resulting equation is equivalent to
 \begin{eqnarray}
   U(\boldsymbol{\beta} )& \doteq & U(\boldsymbol{\beta}, \hat{m}_0(t;\boldsymbol{\beta}) )
       \nonumber \\  [1mm]
 &=& \frac{1}{n} \sum_{i=1}^n \int_0^\tau \pi_i \{ \bZ_i- \bar{\bZ}(t)\} \{\hat{m}_0(t;\boldsymbol{\beta}) d N_i(t)  \nonumber \\
 & & -Y_i(t)\exp(-\boldsymbol{\beta}^{\T} \bZ_i)d t \}=0.  \label{esteq}
 \end{eqnarray}

 The resulting estimator is denoted by
 $\hat{\bbeta}$, which has several good properties such as
 consistency and asymptotic normality.

 In order to derive the large sample properties of $\hat{\bbeta}$, some notations and regularity conditions are needed.
Denote  $H(t|\mathbf{Z})=\Pr (\tilde{T} \geq t |\mathbf{Z})$ and $C_n(t) = n^{-1} \sum_{i=1}^n \pi_i Y_i(t)$. Let
\begin{eqnarray*}
\tilde{\bZ}(t) = \frac{S_n(t)}{C_n(t)} \int_0^t
   \frac{ n^{-1} \sum_{j=1}^n \pi_j [\bZ_j - \bar{\bZ}(u)] d N_j(u)}{ S_n(u) },
\end{eqnarray*}
and let $\mu_{\mathbf{z}}(t)$ and $\tilde{\mu}_{\bz}(t)$ be the limits of $\bar{\mathbf{Z}}(t)$ and $\tilde{\bZ}(t)$,
respectively.

We give the following regularity conditions:
\begin{enumerate}
   \item[C1] $\sup$ supp(F) $\leq$ $\sup$ supp(G), where $F(\cdot)$ and $G(\cdot)$
             are the distribution functions of $T$ and $C$, respectively;
   \item[C2] $\sup_i \| \bZ_i \| < \infty$, where $ \| \bu \|$ denote the Euclidean norm of vector variable $\bu$;
   \item[C3] $m_*(t)$ is continuously differentiable on $[0,\tau]$;
   \item[C4]  $A  = \int_0^\tau \E  \left[  \left\{\mathbf{Z}- \mu_{\mathbf{z}}(t)\right\}^{\otimes 2}
            \exp(-\boldsymbol{\beta}_*^{\T}\textbf{Z}) H(t|\mathbf{Z})\right]d t$ is nonsingular, where $a^{\otimes 2}$ denotes $aa^{\T}$ for a vector $a$.
   \item[C5] $\tilde{n}/n$ converges to $p$ as $n \rightarrow \infty$, where $p_0 \leq p \leq 1$ for some $p_0 > 0$.
 \end{enumerate}

Condition C1 is imposed to ensure that the mean residual life function is estimable, otherwise the MRLF of $T$ may not be estimable
at some points. From the technical arguments, this assumption also saves us from lengthy technical discussion of the tail behavior of the limiting
distributions. Such an assumption has been adopted by other investigators in regression analysis of MRLF, see, for
example, \cite{Chen2005}, \cite{Chen2005b}, \cite{Chen2006} and \cite{sun2009class}. This assumption may not hold if the survival time has an extremely long tail. It may also fail when the follow-up period is too short or when the tail is subject to administrative censoring. However, in well-designed clinical studies
with a nonzero event rate and long follow-up, this assumption is reasonable.

To introduce our results, let
\begin{eqnarray*}
  & & \Sigma  = \Sigma_1+\Sigma_2,   \\ [1mm]
  & & \Sigma_1 = \E \left[\int_0^\tau \left\{\mathbf{Z}_1-\mu_{\mathbf{z}}(t)- \tilde{\mu}_{\bz}(t) \right\}m_*(t)d M_1(t) \right]^{\otimes 2},\\
  & & \Sigma_2 = \frac{1-p}{p} \E \bigg[\int_0^\tau  (1-\delta_1)\left\{\mathbf{Z}_1-\mu_{\mathbf{z}}(t)- \tilde{\mu}_{\bz}(t) \right\} \\
  & & ~~~~~~~~~~~~~ m_*(t) d M_1(t) \bigg]^{\otimes 2}\\
  & & ~~~~~~~~- \frac{1-p}{p} \bigg(\E \bigg[\int_0^\tau  (1-\delta_1)\left\{\mathbf{Z}_1-\mu_{\mathbf{z}}(t)- \tilde{\mu}_{\bz}(t) \right\}\\
  & & ~~~~~~~~~~~~~ m_*(t) d M_1(t) \bigg]     \bigg)^{\otimes 2}.
\end{eqnarray*}

\begin{theorem}
\label{theo1}
Suppose conditions C1--C5 hold, then \\
(i) $\hat{\bbeta}$ and $\hat{m}_0(t)$ always exist and are consistent. \\
(ii) $n^{1/2} ( \hat{\boldsymbol{\beta}}-\boldsymbol{\beta}_* ) $ is asymptotic normal with mean zero and a variance-covariance matrix
$A^{-1} \Sigma (A^{-1})^{\T}$.
Moreover, $A$ and $\Sigma$ can be consistently estimated by
$\hat{A}$ and $\hat{\Sigma}$ respectively, where
\begin{eqnarray*}
& & \hat{A}= \frac{1}{n} \sum_{i=1}^n \pi_i \int_0^\tau \{ \bZ_i - \bar{\bZ} (t) \}^{\otimes 2} Y_i(t)
    \exp( - \hat{\bbeta}^{\T} \bZ_i) d t, \\
& & \hat{\Sigma}= \hat{\Sigma}_1 + \hat{\Sigma}_2,\\
& &  \hat{\Sigma}_1 = \frac{1}{n} \sum_{i=1}^n \pi_i \int_0^\tau \{ \bZ_i - \bar{\bZ} (t) - \tilde{\bZ}(t) \}^{\otimes 2} \\
& & ~~~~~~ Y_i(t) \hat{m}_0(t; \bbeta) \{ \exp( - \hat{\bbeta}^{\T} \bZ_i) d t + d \hat{m}_0(t; \bbeta) \}, \\
& & \hat{\Sigma}_2 =  \frac{1-\hat{p} }{ \hat{p} }  \frac{1}{n} \sum_{i=1}^n \bigg[\int_0^\tau  (1-\delta_i) \frac{ \xi_i }{ \hat{p} }
            \left\{\mathbf{Z}_i - \bar{\bZ}(t) -\tilde{\bZ}(t) \right \} \\
& & ~~~~~~~~~~~~~~\hat{m}_0(t; \bbeta) d \hat{M}_i(t) \bigg]^{\otimes 2}\\
& & ~~~~~ - \frac{1-\hat{p} }{ \hat{p} } \bigg(  \frac{1}{n} \sum_{i=1}^n  \bigg[\int_0^\tau  (1-\delta_i) \frac{ \xi_i }{ \hat{p} }
        \left\{\mathbf{Z}_i - \bar{\bZ}(t) - \tilde{\bZ}(t) \right\} \\
& & ~~~~~~~~~~~~~~ \hat{m}_0(t; \bbeta)  d \hat{M}_i(t) \bigg] \bigg)^{\otimes 2}.
 \end{eqnarray*}
(iii) $n^{1/2} \{ \hat{m}_0(t)-m_*(t) \} (0 \leq t \leq \tau)$ converges weakly to a mean zero Gaussian process with the covariance function that
will be given in the Appendix.
\end{theorem}
The proof of Theorem $\ref{theo1}$ is given in the Appendix.

 Although $\hat{\bbeta}$ has properties such as consistency and asymptotic normality that can be used to make valid inferences about $\bbeta_*$,
the ad hoc nature of $U(\bbeta)$ would not lead to efficient estimators, however.
Note that equation (\ref{esteq}) is constructed via the method of moments, one of the shortcomings for the method
of moments is that it may not necessarily be efficient. To improve the efficiency, we explore the following approaches via two aspects.

One is that by adding proper weight functions.
The following weighted version of the estimating equations can be used to estimate $\bbeta_*$:
\begin{eqnarray}
& & \frac{1}{n} \sum_{i=1}^n \int_0^\tau \pi_i W_i(t) \{ \bZ_i- \bar{\bZ}(t)\} \{\hat{m}_0(t;\boldsymbol{\beta}) d N_i(t)  \nonumber \\
& &~~~~~ - Y_i(t)\exp(-\boldsymbol{\beta}^{\T} \bZ_i)d t \}=0,   \label{weighteq}
\end{eqnarray}
where $W_i(t)$ is a possibly data-dependent and $\mathcal{F}_t-$measurable weight function which converges uniformly to some deterministic function $w(t)$ almost surely.
Denote the solution to the  equation (\ref{weighteq}) as $\hat{\bbeta}_w$,
 by using the technique in the Appendix, $\hat{\bbeta}_w$ is shown to be consistent and asymptotically normal
with asymptotic variance $n^{-1} A_w^{-1} \Sigma_w A_w^{-1}$, where
\begin{eqnarray*}
        A_w &=&  \int_0^\tau \E  \left[ w(t) \left\{\mathbf{Z}- \mu_{\mathbf{z}}(t)\right\}^{\otimes 2}
                 \exp(-\boldsymbol{\beta}_*^{\T}\textbf{Z}) H(t|\mathbf{Z})\right]d t,    \\ [1mm]
\Sigma_w    &=&  \Sigma_{w1} + \Sigma_{w2},   \\ [1mm]
\Sigma_{w1} &=&  \E \bigg[\int_0^\tau w(t) \left\{\mathbf{Z}_1-\mu_{\mathbf{z}}(t)- \tilde{\mu}_{\bz}(t) \right\} \\
            & & ~~~~~ m_*(t)d M_1(t) \bigg]^{\otimes 2},\\
\Sigma_{w2} &=&  \frac{1-p}{p} \E \bigg[\int_0^\tau w(t)  (1-\delta_1)\left\{\mathbf{Z}_1-\mu_{\mathbf{z}}(t)- \tilde{\mu}_{\bz}(t)\right\} \\
            & & ~~~~~ m_*(t)d M_1(t) \bigg]^{\otimes 2}\\
            & &  - \frac{1-p}{p} \bigg(\E \bigg[\int_0^\tau w(t) (1-\delta_1)\left\{\mathbf{Z}_1-\mu_{\mathbf{z}}(t)- \tilde{\mu}_{\bz}(t) \right\}   \\
            & & ~~~~~ m_*(t) d M_1(t) \bigg]     \bigg)^{\otimes 2}.
\end{eqnarray*}
\cite{Chen2005} has used the Cauchy--Schwarz inequality to prove that
\begin{eqnarray}   \label{weight}
W_i(t) = \frac{ \exp(-\bbeta^{\T} \bZ_i) }{ \hat{m}_0(t) \{ \exp(-\bbeta^{\T} \bZ_i) + \hat{m}_0^\prime (t) \} },
\end{eqnarray}
can improve the estimation efficiency. Hence we also suggest this choice of $W_i(t)$, which actually decrease the estimated variance in our simulation results.

The other is through the stratified case-cohort design. For the stratified case-cohort design, the general idea is the same with the classical one. That is, the covariate $\bZ_i$ is observed only when subject $i$ is a failure or from the subcohort. The observed data can be summarized as $\{(\tilde{T}_i, \delta_i, \xi_i, [\delta_i+(1-\delta_i)\xi_i]\bZ_i), i=1, \ldots, n\}$, where $\xi_i$ is the subcohort indicator.
But the sampling scheme to choose the subcohort is no longer simple random sampling.
If an individual characteristics $\bZ^*$ being highly correlated with $\bZ$ is available for all the subjects in the cohort,
\cite{Nan2004} suggested that selecting the subcohort using stratified sampling based on $\bZ^*$ can improve efficiency in hazard regression models.
We expect that a similar result holds for mean residual life models, which have been supported by the simulation studies in next Section.
Many sampling schemes can be designed for selecting a stratified subcohort.
In this paper, we allow $\xi_i$ to depend on  $\bZ_i^*$, which may involve $\tilde{T}_i$, $\bZ_i$ and some external variables correlated with $\tilde{T}_i$ and $\bZ_i$, and the $\xi_i's$ are independent Bernoulli variables with possibly unequal success probabilities.
Let $p_i=\Pr(\xi_i=1)=p(\bZ_i^*)$ be the probability to be chosen in the subcohort, where $p(\bZ_i^*)$ is a function mapping the sample space of $\bZ^*$ to $(p_0, 1)$ for some $p_0 > 0$.
Then the weight under the stratified case-cohort design is defined as $\pi_i^s=\delta_i + (1-\delta_i)\xi_i/p_i$. The resulting estimating equations are
\begin{eqnarray}
& & \sum_{i=1}^n \pi_i^s \big[ m_0(t)dN_i(t)- Y_i(t) \big\{\exp(-\boldsymbol{\beta}^{\T}\textbf{Z}_i) dt  \nonumber \\
& &~~~~~ + d m_0(t) \big\}\big]=0  ~~~(0\leq t\leq \tau),    \label{eqs1}   \\
& & \sum_{i=1}^n \int_0^\tau \pi_i^s \textbf{Z}_i\big[ m_0(t)dN_i(t)- Y_i(t) \big\{\exp(-\boldsymbol{\beta}^{\T}\textbf{Z}_i) dt \nonumber \\
& &~~~~~ + d m_0(t) \big\}\big]=0.   \label{eqs2}
\end{eqnarray}

Note that when there is only a single stratum, the independent Bernoulli sampling proposed for selecting the subcohort in the stratified case-cohort design does not reduce to sampling without replacement.
The size of the subcohort $\tilde{n} = \sum_{i=1}^n \xi_i$ is random, if $\frac{1}{n} \sum_{i=1}^n p_i$ converges to the limiting subcohort proportion $p \in (0, 1]$ in probability, so does $\tilde{n}/n$.
The stratified case-cohort design keeps the independent structure while the classical case-cohort study
does not, which results in the different proofs of their asymptotic properties.
In fact, following the arguments in the Appendix and those in  \cite{Kulich2000},
we can prove that the resulting estimator from (\ref{eqs1}) and (\ref{eqs2}), denote as  $\hat{\bbeta}_s$, is a consistent and asymptotically normal estimator of $\bbeta_*$,
the asymptotic variance of  $\hat{\bbeta}_s$ is $ A^{-1} \Sigma_s (A^{-1})^{\T}$,
where $\Sigma_s  = \Sigma_{s1}+\Sigma_{s2}$ with
\begin{eqnarray*}
   & & \Sigma_{s1} = \E \left[\int_0^\tau \left\{\mathbf{Z}_1-\mu_{\mathbf{z}}(t)- \tilde{\mu}_{\bz}(t) \right\}m_*(t) d M_1(t) \right]^{\otimes 2},\\
  & & \Sigma_{s2} = \frac{1-p}{p} \E \bigg[\int_0^\tau  (1-\delta_1)\left\{\mathbf{Z}_1-\mu_{\mathbf{z}}(t)- \tilde{\mu}_{\bz}(t) \right\} \\
  & & ~~~~~m_*(t) d M_1(t) \bigg]^{\otimes 2}.
\end{eqnarray*}

\section{Simulation Studies}

\renewcommand\arraystretch{1.2}
\begin{table*}[!hbp]
\caption{Simulation results when the censoring rate is approximately 70\%}
\label{tabn1}
\begin{tabular*}{\textwidth}{@{\extracolsep{\fill}}ccccccccccccccc}
\hline\hline
&~~~~~~~~~~~ & ~~~~~~~ & \multicolumn{3}{c}{Full}& ~~~~~~~ &\multicolumn{3}{c}{Subcohort:200} & ~~~~~~~& \multicolumn{3}{c}{Subcohort:300} &  \\
         \cline{4-6}   \cline{8-10}   \cline{12-14}
& & & $\beta_1$ & & $\beta_2$ &  & $\beta_1$ & & $\beta_2$ & & $\beta_1$ & & $\beta_2$ &  \\
\hline
&      &  & \multicolumn{12}{c} {$ m_0(t)= ( -0.5 t + 0.5 )_{+}, \bbeta = (0,0)^{\T} $ } \\
& BIAS &  & 0.005 & & 0.001   & & 0.004 & & 0.007 &   & 0.006 & &-0.001 & \\
& SD   &  & 0.057 & & 0.100   & & 0.097 & & 0.148 &   & 0.083 & & 0.138 & \\
& SE   &  & 0.060 & & 0.104   & & 0.088 & & 0.153 &   & 0.079 & & 0.137 & \\
& CP   &  & 94.8  & & 96.0    & & 91.6  & & 94.4  &   & 93.8  & & 94.2 & \\
& RE   &  & 1.00  & & 1.00    & & 0.35  & & 0.46  &   & 0.48  & & 0.53 & \\
  \hline
&      &  & \multicolumn{12}{c} {$ m_0(t)= ( -0.5 t + 0.5 )_{+}, \bbeta = (0.2,0.2)^{\T} $ } \\
& BIAS &  & 0.002 & &-0.002   & & 0.005 & &-0.003 &   &-0.004 & &-0.004 & \\
& SD   &  & 0.044 & & 0.077   & & 0.067 & & 0.109 &   & 0.060 & & 0.103 & \\
& SE   &  & 0.049 & & 0.086   & & 0.070 & & 0.121 &   & 0.063 & & 0.109 & \\
& CP   &  & 98.0  & & 97.0    & & 96.0  & & 97.6  &   & 97.0  & & 97.2  & \\
& RE   &  & 1.00  & & 1.00    & & 0.43  & & 0.50  &   & 0.58  & & 0.59  & \\
  \hline
&      &  & \multicolumn{12}{c} {$ m_0(t)= ( -0.5 t + 0.5 )_{+}, \bbeta = (0.5,-0.5)^{\T} $ } \\
& BIAS &  & 0.002 & & 0.007   & & 0.002 & &-0.002 &   & 0.007 & & 0.002 & \\
& SD   &  & 0.058 & & 0.104   & & 0.086 & & 0.160 &   & 0.079 & & 0.142 & \\
& SE   &  & 0.062 & & 0.106   & & 0.093 & & 0.159 &   & 0.083 & & 0.143 & \\
& CP   &  & 97.2  & & 95.8    & & 96.4  & & 95.4  &   & 97.6  & & 95.2  & \\
& RE   &  & 1.00  & & 1.00    & & 0.43  & & 0.43  &   & 0.52  & & 0.54  & \\
  \hline
&      &  & \multicolumn{12}{c} {$ m_0(t)= ( -0.5 t + 1 )_{+}, \bbeta = (0,0)^{\T} $ } \\
& BIAS &  & 0.001 & &-0.003   & &-0.001 & &-0.014 &   & 0.005 & & 0.003 & \\
& SD   &  & 0.058 & & 0.097   & & 0.090 & & 0.163 &   & 0.077 & & 0.141 & \\
& SE   &  & 0.060 & & 0.104   & & 0.089 & & 0.154 &   & 0.080 & & 0.138 & \\
& CP   &  & 94.8  & & 96.6    & & 93.6  & & 92.6  &   & 96.0  & & 93.0  & \\
& RE   &  & 1.000 & & 1.000   & & 0.416 & & 0.351 &   & 0.574 & & 0.470 & \\
  \hline
&      &  & \multicolumn{12}{c} {   $ m_0(t)= ( -0.5 t + 1 )_{+}, \bbeta = (0.2, 0.2)^{\T} $   } \\
& BIAS &  & -0.003  & & -0.001    & & -0.008  & &  0.007  &   & -0.003  & &  0.003  & \\
& SD   &  &  0.044  & &  0.082    & &  0.063  & &  0.121  &   &  0.062  & &  0.106  & \\
& SE   &  &  0.050  & &  0.086    & &  0.070  & &  0.120  &   &  0.063  & &  0.109  & \\
& CP   &  &  97.4   & &  96.0     & &  97.2   & &  94.4   &   &  95.6   & &  95.8   & \\
& RE   &  &  1.00   & &  1.00     & &  0.49   & &  0.46   &   &  0.50   & &  0.60   & \\
  \hline
&      &  & \multicolumn{12}{c} {$ m_0(t)= ( -0.5 t + 1 )_{+}, \bbeta = (0.5, -0.5)^{\T} $ } \\
& BIAS &  & 0.002 & & 0.019   & & 0.003 & & 0.027 &   & 0.002 & & 0.020 & \\
& SD   &  & 0.058 & & 0.102   & & 0.088 & & 0.156 &   & 0.080 & & 0.147 & \\
& SE   &  & 0.062 & & 0.106   & & 0.093 & & 0.159 &   & 0.083 & & 0.142 & \\
& CP   &  & 95.6  & & 95.8    & & 95.6  & & 95.2  &   & 94.6  & & 93.0  & \\
& RE   &  & 1.00  & & 1.00    & & 0.43  & & 0.42  &   & 0.51  & & 0.47  & \\
  \hline
\end{tabular*}
\begin{tablenotes}
\item BIAS, the empirical bias; SD, the empirical standard deviation; SE, the mean of
estimated standard error; CP, the empirical coverage probability of $95\%$ confidence interval;
RE, the empirical relative efficiencies, calculated by the ratio of sample variance with the full estimators as a reference.
\end{tablenotes}
\end{table*}

\renewcommand\arraystretch{1.2}
\begin{table*}[!hbp]
\caption{Simulation results when the censoring rate is approximately 80\%}
\label{tabn2}
\begin{tabular*}{\textwidth}{@{\extracolsep{\fill}}ccccccccccccccc}
\hline\hline
&~~~~~~~~~~~ & ~~~~~~~ & \multicolumn{3}{c}{Full}& ~~~~~~~ &\multicolumn{3}{c}{Subcohort:200} & ~~~~~~~& \multicolumn{3}{c}{Subcohort:300} &  \\
         \cline{4-6}   \cline{8-10}   \cline{12-14}
& & & $\beta_1$ & & $\beta_2$ &  & $\beta_1$ & & $\beta_2$ & & $\beta_1$ & & $\beta_2$ &  \\
\hline
&      &  & \multicolumn{12}{c} {$ m_0(t)= ( -0.5 t + 0.5 )_{+}, \bbeta = (0,0)^{\T} $ } \\
& BIAS &  & 0.003 & & 0.003   & &-0.003 & & 0.010 &   & 0.006 & &-0.012 & \\
& SD   &  & 0.066 & & 0.118   & & 0.104 & & 0.186 &   & 0.093 & & 0.171 & \\
& SE   &  & 0.073 & & 0.127   & & 0.105 & & 0.182 &   & 0.095 & & 0.165 & \\
& CP   &  & 97.2  & & 96.8    & & 94.8  & & 93.0  &   & 96.8  & & 94.6  & \\
& RE   &  & 1.00  & & 1.00    & & 0.40  & & 0.40  &   & 0.50  & & 0.48  & \\
  \hline
&      &  & \multicolumn{12}{c} {$ m_0(t)= ( -0.5 t + 0.5 )_{+}, \bbeta = (0.2,0.2)^{\T} $ } \\
& BIAS &  &-0.004 & &-0.004   & &-0.007 & & 0.003 &   &-0.005 & &-0.008 & \\
& SD   &  & 0.054 & & 0.098   & & 0.079 & & 0.139 &   & 0.071 & & 0.129 & \\
& SE   &  & 0.059 & & 0.102   & & 0.083 & & 0.143 &   & 0.075 & & 0.129 & \\
& CP   &  & 97.2  & & 96.0    & & 95.8  & & 95.8  &   & 95.6  & & 93.6  & \\
& RE   &  & 1.00  & & 1.00    & & 0.47  & & 0.50  &   & 0.58  & & 0.57  & \\
  \hline
&      &  & \multicolumn{12}{c} {$ m_0(t)= ( -0.5 t + 0.5 )_{+}, \bbeta = (0.5,-0.5)^{\T} $ } \\
& BIAS &  &-0.002 & & 0.010   & &-0.011 & & 0.015 &   &-0.006 & & 0.012 & \\
& SD   &  & 0.069 & & 0.126   & & 0.107 & & 0.196 &   & 0.093 & & 0.174 & \\
& SE   &  & 0.074 & & 0.128   & & 0.110 & & 0.190 &   & 0.099 & & 0.172 & \\
& CP   &  & 96.2  & & 94.8    & & 93.8  & & 93.0  &   & 95.2  & & 94.0  & \\
& RE   &  & 1.00  & & 1.00    & & 0.41  & & 0.41  &   & 0.55  & & 0.52  & \\
  \hline
&      &  & \multicolumn{12}{c} {$ m_0(t)= ( -0.5 t + 1 )_{+}, \bbeta = (0,0)^{\T} $ } \\
& BIAS &  & 0.002 & &-0.005   & &-0.002 & & 0.001 &   & 0.002 & &-0.007 & \\
& SD   &  & 0.069 & & 0.115   & & 0.111 & & 0.179 &   & 0.094 & & 0.158 & \\
& SE   &  & 0.073 & & 0.127   & & 0.106 & & 0.183 &   & 0.096 & & 0.166 & \\
& CP   &  & 96.4  & & 97.4    & & 92.0  & & 92.8  &   & 93.8  & & 95.8  & \\
& RE   &  & 1.00  & & 1.00    & & 0.38  & & 0.41  &   & 0.54  & & 0.53  & \\
  \hline
&      &  & \multicolumn{12}{c} {$ m_0(t)= ( -0.5 t + 1 )_{+}, \bbeta = (0.2, 0.2)^{\T} $ } \\
& BIAS &  &-0.001 & &-0.010   & &-0.004 & &-0.013 &   &-0.003 & &-0.003 & \\
& SD   &  & 0.056 & & 0.096   & & 0.085 & & 0.147 &   & 0.072 & & 0.123 & \\
& SE   &  & 0.059 & & 0.102   & & 0.082 & & 0.142 &   & 0.075 & & 0.129 & \\
& CP   &  & 96.8  & & 96.0    & & 94.0  & & 93.0  &   & 94.4  & & 95.6  & \\
& RE   &  & 1.00  & & 1.00    & & 0.43  & & 0.43  &   & 0.58  & & 0.59  & \\
  \hline
&      &  & \multicolumn{12}{c} {$ m_0(t)= ( -0.5 t + 1 )_{+}, \bbeta = (0.5, -0.5)^{\T} $ } \\
& BIAS &  &-0.003 & & 0.001   & &-0.003 & & 0.004 &   &-0.009 & & 0.002 & \\
& SD   &  & 0.075 & & 0.133   & & 0.110 & & 0.195 &   & 0.095 & & 0.181 & \\
& SE   &  & 0.074 & & 0.129   & & 0.111 & & 0.191 &   & 0.099 & & 0.172 & \\
& CP   &  & 94.8  & & 93.8    & & 92.6  & & 91.0  &   & 94.4  & & 92.8  & \\
& RE   &  & 1.00  & & 1.00    & & 0.46  & & 0.47  &   & 0.64  & & 0.54  & \\
  \hline
\end{tabular*}
\begin{tablenotes}
\item BIAS, the empirical bias; SD, the empirical standard deviation; SE, the mean of
estimated standard error; CP, the empirical coverage probability of $95\%$ confidence interval;
RE, the empirical relative efficiencies, calculated by the ratio of sample variance with the full estimators as a reference.
\end{tablenotes}
\end{table*}

In this Section, we conduct simulation studies to examine the finite sample properties
of the proposed estimator.

In the first scenario of simulation studies, the event time $T$ is generated from the following proportional mean residual life regression model
\begin{eqnarray*}
m(t|\bZ)=m_0(t) \exp(\beta_1 Z_1 + \beta_2 Z_2),
\end{eqnarray*}
where the covariate $Z_1$ is a Bernoulli random variable with success probability
$0.5$, $Z_2$ is a U(0,1) variable, the true regression parameters $(\beta_{1*}, \beta_{2*})^{\T}$
is respectively set to be $(0,0)^{\T}$, $(0.2, 0.2)^{\T}$ or $(0.5, -0.5)^{\T}$
and the baseline function $m_0(t)$ is taken from
the Hall-Wellner family, in other words,
$m_0(t)=(D_1 t + D_2)^{+}$, where
$D_1 > -1, ~ D_2 > 0,$ and $d^{+}$ denotes $d \cdot I(d \geq 0)$
for any quantity $d$.
Here we consider two scenarios for the baseline function $m_0(t)$.
One is that $D_1=-0.5$
and $D_2=0.5$, the other is under $D_1=-0.5$
and $D_2=1$.
The censoring time $C$ is generated from $\mbox{Exp}(\lambda)$, where $\lambda$ is used to control the censoring proportion.
We set the censoring rate to be $70\%$ or $80\%$ to mimic the low event rate where the case-cohort designs are more likely to be applied.
Care is taken to ensure that the support of $C$ is larger than the support of $T$ for all $\bZ$.

$500$ replications of full cohort data are generated  with the sample size $n=1000$.
For each replication, subcohorts of size $200$ and size $300$
are drawn by simple random sampling, respectively.
The empirical biases (Bias), empirical standard deviations (SD),
average robust standard errors (SE), coverage probabilities of the
$95\%$ confidence intervals (CP) and the empirical relative efficiency (RE)
of the proposed $\hat{\bbeta}$ are reported in the study.
We also report the estimation and inference results based on full cohort
for comparison.
The simulation results are summarized in  Tables $\ref{tabn1}$ and $\ref{tabn2}$, when the censoring rates are approximately $70\%$ and $80\%$, respectively.

It can be seen from the simulation results in Tables $\ref{tabn1}$ and $\ref{tabn2}$ that
the proposed estimates are all essentially
unbiased under two different subcohorts.
The means of estimated standard errors match the empirical standard errors
quite well and the $95\%$ confidence intervals
have reasonable coverage rates.
Compared to the full cohort analysis, the case-cohort designs are less
efficient in estimating the regression coefficients but
the efficiency loss appears to be small.
We also find that the empirical relative efficiencies increase when the size of the subcohort increases.

It would be informative to compare the proposed case-cohort estimator with the estimator based on a simple random sample of the same size as the case-cohort sample. For the sample size $n$, the subcohort size $n_1$, and the censoring rate $c$, the effective case-cohort sample size is $n_1 + (n-n_1)\cdot(1-c)$. When $n=1000, c=0.8$, subcohorts of size 200 and 300 correspond to the effective size 360 and 440, respectively. For the same data generation mechanism as in Table \ref{tabn2}, we present the simulation comparison results in Table \ref{tabr3}. We can see from Table \ref{tabr3} that the relative efficiencies are all larger than 1, showing the advantage of using a case-cohort design instead of a random sampling design. Compared to the same size of random sampling design, the estimated standard errors and the empirical deviations of the case-cohort sample are much closer to each other.

\renewcommand\arraystretch{1.2}
\begin{table*}[!hbp]
\caption{Simulation comparison when the censoring rate is approximately 80\%}
\label{tabr3}
\begin{tabular*}{\textwidth}{@{\extracolsep{\fill}}ccccccccccccccc}
\hline\hline
 &  & \multicolumn{2}{c}{$n$=1000, subcohort:200}&   &\multicolumn{3}{c}{$n$=360} &  &  \multicolumn{2}{c}{$n$=1000, subcohort:300}&   &\multicolumn{3}{c}{$n$=440}  \\
         \cline{3-4}   \cline{6-8}   \cline{10-11}  \cline{13-15}
 &  & $\beta_1$  & $\beta_2$ &  & $\beta_1$ & & $\beta_2$  &  & $\beta_1$  & $\beta_2$ &  & $\beta_1$ & & $\beta_2$  \\
\hline
      &  & \multicolumn{13}{c} {$ m_0(t)= ( -0.5 t + 0.5 )_{+}, \bbeta = (0,0)^{\T} $ } \\
 BIAS &  & 0.001  &-0.011   & &-0.002 & & 0.014 &   & 0.002  &-0.008 & & 0.002 & &-0.001  \\
 SD   &  & 0.110  & 0.181   & & 0.117 & & 0.212 &   & 0.102  & 0.166 & & 0.106 & & 0.193  \\
 SE   &  & 0.105  & 0.182   & & 0.122 & & 0.211 &   & 0.096  & 0.166 & & 0.110 & & 0.191  \\
 CP   &  & 92.2   & 92.6    & & 95.4  & & 94.0  &   & 92.6   & 95.0  & & 95.6  & & 96.0   \\
 RE   &  & 1.33   & 1.35    & & 1.00  & & 1.00  &   & 1.33   & 1.33  & & 1.00  & & 1.00   \\
  \hline
      &  & \multicolumn{13}{c} {$ m_0(t)= ( -0.5 t + 0.5 )_{+}, \bbeta = (0.2,0.2)^{\T} $ } \\
 BIAS &  &-0.001  & 0.005   & &-0.007 & &-0.001 &   &-0.003  & 0.002 & &-0.003 & &-0.009  \\
 SD   &  & 0.082  & 0.140   & & 0.092 & & 0.157 &   & 0.071  & 0.133 & & 0.082 & & 0.145  \\
 SE   &  & 0.082  & 0.142   & & 0.099 & & 0.172 &   & 0.075  & 0.130 & & 0.090 & & 0.156  \\
 CP   &  & 94.8   & 95.4    & & 95.0  & & 95.0  &   & 95.2   & 94.6  & & 97.2  & & 95.8  \\
 RE   &  & 1.45   & 1.45    & & 1.00  & & 1.00  &   & 1.44   & 1.44  & & 1.00  & & 1.00  \\
  \hline
      &  & \multicolumn{13}{c} {$ m_0(t)= ( -0.5 t + 0.5 )_{+}, \bbeta = (0.5,-0.5)^{\T} $ } \\
 BIAS &  &-0.007  & 0.010   & &-0.022 & & 0.018 &   &-0.001  & 0.019 & &-0.012 & & 0.010  \\
 SD   &  & 0.107  & 0.194   & & 0.116 & & 0.205 &   & 0.102  & 0.179 & & 0.113 & & 0.191  \\
 SE   &  & 0.109  & 0.190   & & 0.122 & & 0.212 &   & 0.098  & 0.171 & & 0.111 & & 0.193  \\
 CP   &  & 95.0   & 92.4    & & 94.6  & & 94.2  &   & 94.2   & 94.6  & & 93.8  & & 95.0  \\
 RE   &  & 1.26   & 1.25    & & 1.00  & & 1.00  &   & 1.29   & 1.28  & & 1.00  & & 1.00  \\
  \hline
      &  & \multicolumn{13}{c} {$ m_0(t)= ( -0.5 t + 1 )_{+}, \bbeta = (0,0)^{\T} $ } \\
 BIAS &  &-0.001  &-0.001   & &-0.001 & &-0.013 &   &-0.001  & 0.001 & & 0.002 & &-0.005  \\
 SD   &  & 0.105  & 0.194   & & 0.112 & & 0.207 &   & 0.097  & 0.173 & & 0.106 & & 0.187  \\
 SE   &  & 0.105  & 0.181   & & 0.122 & & 0.212 &   & 0.095  & 0.164 & & 0.111 & & 0.192  \\
 CP   &  & 93.4   & 93.0    & & 96.2  & & 94.2  &   & 93.8   & 93.6  & & 96.2  & & 95.4  \\
 RE   &  & 1.36   & 1.37    & & 1.00  & & 1.00  &   & 1.37   & 1.38  & & 1.00  & & 1.00  \\
  \hline
      &  & \multicolumn{13}{c} {$ m_0(t)= ( -0.5 t + 1 )_{+}, \bbeta = (0.2, 0.2)^{\T} $ } \\
 BIAS &  & 0.002  & 0.001   & & 0.003 & &-0.008 &   & 0.001  &-0.004 & &-0.005 & & 0.005  \\
 SD   &  & 0.083  & 0.139   & & 0.097 & & 0.156 &   & 0.074  & 0.124 & & 0.085 & & 0.147  \\
 SE   &  & 0.083  & 0.143   & & 0.099 & & 0.172 &   & 0.075  & 0.130 & & 0.089 & & 0.154  \\
 CP   &  & 95.2   & 93.4    & & 94.6  & & 95.4  &   & 93.4   & 96.2  & & 95.2  & & 93.8  \\
 RE   &  & 1.45   & 1.45    & & 1.00  & & 1.00  &   & 1.42   & 1.42  & & 1.00  & & 1.00  \\
  \hline
      &  & \multicolumn{13}{c} {$ m_0(t)= ( -0.5 t + 1 )_{+}, \bbeta = (0.5, -0.5)^{\T} $ } \\
 BIAS &  &-0.014  & 0.017   & &-0.018 & & 0.019 &   &-0.004  & 0.017 & &-0.008 & & 0.027  \\
 SD   &  & 0.112  & 0.192   & & 0.116 & & 0.215 &   & 0.101  & 0.176 & & 0.113 & & 0.178  \\
 SE   &  & 0.109  & 0.190   & & 0.123 & & 0.214 &   & 0.098  & 0.170 & & 0.112 & & 0.194  \\
 CP   &  & 92.0   & 92.4    & & 96.0  & & 94.0  &   & 94.4   & 94.8  & & 93.6  & & 96.0  \\
 RE   &  & 1.27   & 1.27    & & 1.00  & & 1.00  &   & 1.31   & 1.30  & & 1.00  & & 1.00  \\
  \hline
\end{tabular*}
\begin{tablenotes}
\item BIAS, the empirical bias; SD, the empirical standard deviation; SE, the mean of
estimated standard error; CP, the empirical coverage probability of $95\%$ confidence interval;
RE, the relative efficiencies, calculated by the ratio of estimated variance with the full estimators as a reference.
\end{tablenotes}
\end{table*}

\renewcommand\arraystretch{1.2}
\begin{table*}[!hbp]
\caption{Comparison of classical and stratified case-cohort design }
\label{tabn5}
\begin{tabular*}{\textwidth}{@{\extracolsep{\fill}}ccccccccccc}
\hline\hline
& &   & \multicolumn{3}{c}{$p_z = 0.5$} & &  \multicolumn{3}{c}{$p_z = 0.3$} &  \\
         \cline{4-6}   \cline{8-10}
& &  & Full  &SRS  & STRAT &  & Full  &SRS  & STRAT &    \\
\hline
&      &  & \multicolumn{8}{c} {$ m_0(t)= 1, ~~\beta_0 = 0 $ } \\
& BIAS &  & 0.001    & 0.004  &-0.001 & &-0.009    &-0.011  &-0.015 &\\
& SD   &  & 0.094    & 0.110  & 0.094 & & 0.103    & 0.125  & 0.127 &\\
& SE   &  & 0.098    & 0.115  & 0.098 & & 0.107    & 0.125  & 0.123 &\\
& CP   &  & 96.0     & 96.4   & 96.0  & & 95.2     & 95.4   & 94.4  &\\
\hline
&      &  & \multicolumn{8}{c} {$ m_0(t)= 1, ~~\beta_0 = 0.2 $ } \\
& BIAS &  & 0.002    &-0.002  &-0.003 & & 0.013    & 0.010  &-0.001 &\\
& SD   &  & 0.099    & 0.113  & 0.122 & & 0.114    & 0.125  & 0.130 &\\
& SE   &  & 0.100    & 0.118  & 0.117 & & 0.118    & 0.137  & 0.134 &\\
& CP   &  & 95.0     & 96.0   & 94.8  & & 96.0     & 96.6   & 95.4 &\\
\hline
&      &  & \multicolumn{8}{c} {$ m_0(t)= 1, ~~\beta_0 = 0.5 $ } \\
& BIAS &  & 0.003    &-0.005  &-0.007 & & 0.051    & 0.034  & 0.031 &\\
& SD   &  & 0.110    & 0.123  & 0.132 & & 0.142    & 0.142  & 0.151 &\\
& SE   &  & 0.107    & 0.128  & 0.128 & & 0.141    & 0.162  & 0.161 &\\
& CP   &  & 94.8     & 94.4   & 93.8  & & 96.0     & 96.2   & 96.4 &\\
\hline
&      &  & \multicolumn{8}{c} {$ m_0(t)= 1+t, ~~\beta_0 = 0 $ } \\
& BIAS &  & 0.005    & 0.016  & 0.013 & &-0.006    &-0.006  &-0.017 &\\
& SD   &  & 0.173    & 0.208  & 0.206 & & 0.191    & 0.225  & 0.229 &\\
& SE   &  & 0.170    & 0.200  & 0.202 & & 0.184    & 0.218  & 0.213 &\\
& CP   &  & 93.6     & 93.0   & 93.4  & & 94.2     & 95.0   & 92.8 &\\
\hline
&      &  & \multicolumn{8}{c} {$ m_0(t)= 1+t, ~~\beta_0 = 0.2 $ } \\
& BIAS &  & -0.021    &-0.025  &-0.030 & &-0.014    &-0.014  &-0.022 &\\
& SD   &  & 0.170     & 0.206  & 0.209 & & 0.192    & 0.222  & 0.218 &\\
& SE   &  & 0.169     & 0.200  & 0.201 & & 0.199    & 0.235  & 0.229 &\\
& CP   &  & 94.0      & 94.6   & 92.4  & & 96.0     & 96.2   & 95.4 &\\
\hline
&      &  & \multicolumn{8}{c} {$ m_0(t)= 1+t, ~~\beta_0 = 0.5 $ } \\
& BIAS &  &-0.070    &-0.081  &-0.078 & &-0.027    &-0.046  &-0.053 &\\
& SD   &  & 0.181    & 0.199  & 0.214 & & 0.212    & 0.253  & 0.241 &\\
& SE   &  & 0.178    & 0.202  & 0.210 & & 0.218    & 0.255  & 0.252 &\\
& CP   &  & 92.4     & 93.2   & 93.6  & & 94.8     & 93.4   & 93.4  &\\
  \hline
\end{tabular*}
\begin{tablenotes}
\item SRS, the simple random sampling;
\item STRAT, stratified case-cohort design with $\eta = \nu = 0.7$.
\end{tablenotes}
\end{table*}

In the second scenario of simulation studies, we explore two efficiency improving approaches.
One method is to use the weighted function (\ref{weight}), the other method is to use the stratified case-cohort design.
Here the event time $T$ is generated from
\begin{eqnarray*}
m(t| Z)=m_0(t) \exp( \beta_0 Z ),
\end{eqnarray*}
where the covariate $Z = 2 \cdot \mbox{Bernoulli}(p_z)-1$ with $p_z=0.3$ or $0.5$, the true regression parameter $\beta_0$ is set to be $0$, $0.2$ or $0.5$,
the baseline function $m_0(t)$ is taken from $m_0(t)=1$ or $m_0(t)=1+t$, respectively.
$C$ is generated as described in the first scenario of simulations, but the censoring rates are set approximately 90\% this time.
The sample size and the subcohort size are almost equal to the first ones.
The simulation results using the weighted estimating equations are omitted here because the weighted estimators don't show  significant efficiency improvement in many cases.
Meanwhile the weighted case-cohort estimators are still not efficient as the full estimators, this is in accordance with our expectation since we only use the uncensored data and the subcohort data under the case-cohort design. Similar phenomenon has also been founded in \cite{Lu2006}.

In the stratified case-cohort design simulation, we define the distribution of $Z^*$ by $\eta=\Pr(Z^*=1|Z=1)$ and $\nu=\Pr(Z^*=-1|Z=-1)$, where
$(\eta, \nu)$ is chosen as (0.7, 0.7). Thus  $Z^*=2 \cdot \mbox{Bernoulli}( (1-\nu)(1-p_z)+\eta \cdot p_z )-1$.
The subcohort is a stratified sample selected by independent Bernoulli sampling with selection probability $p(Z^*)$ chosen so that
approximately equal numbers of subjects are selected from the two strata, i.e. $\{Z^*=1\}$ and $\{Z^*=-1\}$.
Simulation results comparing the full, classical, stratified cohort are given in Table \ref{tabn5}.
In general, the stratified case-cohort design behaves better than the classical one when the correlation between $Z$ and $Z^*$ exists.
But when $Z$ and $Z^*$ are uncorrelated, the classical case-cohort design will do slightly better than the stratified one.

\section{A Real Data Example}

\begin{table*}[!hbp ]
\caption{Regression Analyses of Time from the First Employment to the Nasal Sinus Cancer
Death for the South Welsh Nickel Refiners Study}
\label{table:2}
\begin{tabular*}{\textwidth}{@{\extracolsep{\fill}}ccccccccc}
\hline\hline
   Parameter         &  & Full Cohort &  &Case-cohort& & Multiple Case-cohort &    \\
\hline
  $\log$ (AFE-10)    &  &             &  &           & &       &      \\
       Est.          &  &  $-$0.096     &  &  $-$0.060   & &$-$0.082 &    \\
       S.E.          &  &   0.007     &  &   0.014   & & 0.016 &     \\
      P Value        &  &  $<$0.0001  &  &$<$0.0001  & &$<$0.0001&    \\
  (YFE-1915)/10      &  &             &  &           & &       &         \\
       Est.          &  &  $-$0.009     &  & $-$0.024    & &$-$0.007 &     \\
       S.E.          &  &  0.013      &  &  0.019    & & 0.018 &    \\
      P Value        &  &   0.475     &  &  0.195    & & 0.691 &    \\
  (YFE-1915)$^2$/100 &  &             &  &           & &       &      \\
       Est.          &  &  0.090      &  &  0.064    & & 0.072 &     \\
       S.E.          &  &  0.026      &  &  0.035    & & 0.036 &     \\
      P Value        &  &  0.0005     &  &  0.066    & & 0.045 &    \\
  $\log$(EXP+1)      &  &             &  &           & &       &       \\
       Est.          &  &  $-$0.057     &  & $-$0.057    & &$-$0.050 &     \\
       S.E.          &  &   0.013     &  &  0.016    & & 0.018 &     \\
      P Value        &  &  $<$0.0001  &  &  0.0005   & & 0.006 &     \\
\hline
\end{tabular*}
\begin{tablenotes}
\item Est., Parameter Estimate; S.E., Standard Error.
\end{tablenotes}
\end{table*}

In this Section, we apply the proposed case-cohort analysis approach under the
proportional mean residual life model to the South Welsh nickel refiners study.
In this study, men employed in a nickel refinery in South Welsh were
investigated to determine the risk of developing carcinoma of the bronchi
and nasal sinuses which is associated with the refining of nickel.
The cohort was identified using the weekly payrolls of the company and followed
from the year 1934 until 1981.
The complete records of 679 workers employed before 1925
can be obtained from the Appendix \Rmnum{8}
in \cite{Breslow1987}.
Among the full cohort, there were 56 deaths from
cancer of the nasal sinus until 1981.
The event rate for this study is quite low
and hence the case-cohort design is more likely to be applied.

\cite{Breslow1987} used the Cox model to analyse
the mortality data on nasal sinus cancer.
They considered the survival time to be years since first employment and found three significant risk factors:
AFE (age at first employment), YFE (year at fist employment) and
EXP (exposure level).
\cite{Lin1993} fitted the same model to the data
obtained from a ``hypothetical" case-cohort design which was randomly selected 100 subcohort members from the entire cohort.
In this paper, we  fit respectively the mean residual life model
to the full cohort and the case-cohort in which a subcohort with size 100 is drawn by simple random sampling.
To avoid the random effects, we repeat case-cohort analysis 100 times (choose the subcohort 100 times), and name this as multiple case-cohort analysis.
The covariates transformations adopted by \cite{Breslow1987}
are reserved here.
Specifically, we consider four covariates:
$\log$(AFE-10), (YFE-1915)/10, (YFE-1915)$^2$/100 and $\log$(EXP+1).
In Table $\ref{table:2}$, we present estimates, standard errors, and $p$-values of the regression coefficients
under proportional mean residual life model for the full-cohort and case-cohort analysis. The average results for multiple case-cohort analysis are also shown in Table \ref{table:2}.

In general, the proposed case-cohort estimates for each covariate are close
to the corresponding full-cohort estimates, the multiple case-cohort estimates are closer.
This implies that the case-cohort analysis results are convincing.
Both the full cohort and
case-cohort analysis results under proportional mean residual life model
indicate that the covariates $\log$(AFE-10), (YFE-1915)$^2$/100 and $\log$(EXP+1)
 have significant influence on the survival time,
which in accordance with the results indicated by the Cox model under full cohort.
The results of the estimated coefficients show two different scenes: the mean residual life decreases with the individual's AFE or EXP increasing, but has performance for (YFE-1915)$^2$/100 on the contrary.
It is interesting that the regression coefficients for those significant covariates based
on the proportional mean residual life model have opposite signs to their counterparts under the Cox model.
This phenomenon also appeared in the real data analysis of \cite{Chen2005} because of the different
link between the proportional mean residual life model and the Cox model.

\section{Concluding Remarks}

In this paper, we propose some new estimating functions to deal with case-cohort data under proportional mean residual life model.
Appropriate weighted availability indicators are defined when the subcohort is drawn by simple random sampling. The large sample properties of the proposed
estimators are established.

For case-cohort analysis under the Cox model, \cite{Kulich2004} suggested estimating the
 sampling probability $p$ with a weight estimator to achieve further efficiency. Inspired by the idea,
we can consider the weighted estimator proposed in \cite{Chen2009} for $p$ by
\begin{eqnarray*}
 \hat{p}(t)=\frac{ \sum_{i=1}^n \xi_i (1-\delta_i) c_i(t)}{ \sum_{i=1}^n (1-\delta_i) c_i(t)},
\end{eqnarray*}
 where $c_i(t)$ is possibly time-dependent and satisfy some regularity conditions.
Under the Cox model, various versions of the weight $c_i(t)$ for estimating the sampling probability $p$, including both time-constant and
 time-dependent weights, has been suggested by \cite{Chen1999}, \cite{Borgan2000} and \cite{Kulich2004}. Two common choices for $c_i(t)$ are $c_i(t)=1$
 and $c_i(t)=Y_i(t)$, where $Y_i(t)$ is defined in Section 2. Similarly, in mean residual life model, a Horvitz-Thompson type weighted function $\pi_i(t)= \delta_i + ( 1-\delta_i )\xi_i / \hat{p}(t)$
can be considered to replace $\pi_i$ in the estimating equations (\ref{eq113}) and (\ref{eq002}).

\section{Appendix}

{\it Proof of the Theorem $\ref{theo1}$:} \\
(i) Note that
\begin{eqnarray*}
& & \hat{m}_0(t; \bbeta) \\ [1mm]
&=& \frac{1}{S_n(t)} \int_t^\tau S_n(u) \frac{ \sum_{i=1}^n \pi_i Y_i(u)  \exp(-\bbeta^{\T} \bZ_i)  }
                                { \sum_{i=1}^n \pi_i Y_i(u) } d u  \\ [1mm]
&=& \frac{1}{ E[S(t|\bZ)] } \int_t^\tau E[S(u|\bZ)]  \frac{ \int_\bz S(u|\bz) \exp(-\bbeta^{\T} \bz) d F_\bz (\bz)  }{E[S(u|\bZ)]} du \\ [1mm]
& &   + o_p(1)    \\ [1mm]
&=& \frac{ E[ \exp(-\bbeta^{\T} \bZ) \int_t^\tau S(u|\bZ) d u ] }{ E[S(t|\bZ)] }  + o_p(1)    \\ [1mm]
& \doteq & m_0(t; \bbeta) + o_p(1),
\end{eqnarray*}
where $F_\mathbf{z}(\mathbf{z}) $ is the distribution function of $\mathbf{Z}$ and
$S(t|\bZ)$ is the survival function of $T$ given $\bZ$.
This implies that $\hat{m}_0(t; \bbeta)$ converges to $m_0(t; \bbeta)$ uniformly in $t \in [0, \tau]$ and $\bbeta$ in a compact
set which contains the true parameter $\bbeta_*$,
and $m_0(t; \bbeta_*) = m_*(t)$.
Therefore, to prove the existence of $\hat{\bbeta}$ and $\hat{m}_0(t)$, it suffices to show that there exists a solution to $U(\bbeta) = 0$.
By differentiating $ U( \boldsymbol{\beta} ) $ with respect to $\boldsymbol{\beta}$, we have
\begin{eqnarray*}
& & \hat{A}(\bbeta_*)   \\ [1mm]
& \doteq &  \left.\frac{\partial U(\boldsymbol{\beta} )} {\partial \boldsymbol{\beta}} \right|_{\boldsymbol{\beta}=\boldsymbol{\beta}_*}
                            \\ [1mm]
&=& \frac{1}{n} \sum_{i=1}^n \int_0^\tau  \pi_i  \left\{\mathbf{Z}_i-\bar{\mathbf{Z}}(t) \right\}
     \big\{ -m_*(t) \mu_{\mathbf{z}}(t)  d N_i(t) \\ [1mm]
& &~~~~~ +\mathbf{Z}_i \exp(-\boldsymbol{\beta}_*^{\T}\textbf{Z}_i) Y_i(t) d t \big\}^{\T} +o_p(1)    \\ [1mm]
&=& \frac{1}{n} \sum_{i=1}^n \int_0^\tau \pi_i \left\{\mathbf{Z}_i-\bar{\mathbf{Z}}(t) \right\}
     \left\{ -m_*(t) \mu_{\mathbf{z}}(t)\right\}^{\T} d M_i(t) \\[1mm]
& & +\frac{1}{n} \sum_{i=1}^n \int_0^\tau  \pi_i   \left\{\mathbf{Z}_i-\bar{\mathbf{Z}}(t) \right\}
     \left\{\mathbf{Z}_i-\mu_{\mathbf{z}}(t)\right\}^{\T}  \\ [1mm]
& & ~~~~~ \exp(-\boldsymbol{\beta}_*^{\T}\textbf{Z}_i) Y_i(t) d t +o_p(1)\\ [1mm]
&=& A+o_p(1),
\end{eqnarray*}
this implies that $\hat{A}(\bbeta_*)$ converges in probability to a nonrandom matrix $A$. Since $U(\bbeta_*) \rightarrow 0$ almost surely, and $A$ is nonsingular
by the regularity condition C4, the convergence of $\hat{A}(\bbeta_*)$ implies that we can find a small neighborhood of $\bbeta_*$ in which $\hat{A}(\bbeta_*)$
is nonsingular when $n$ is large enough. Hence it follows from the inverse function theorem that within a small neighborhood of $\bbeta_*$,
there exists a solution $\hat{\bbeta}$ to $U(\hat{\bbeta}) = 0$ for sufficiently large $n$.
Notice that $\hat{\bbeta}$ is strongly consistent to $\bbeta_*$, then it follows from the uniform convergence of $\hat{m}_0(t; \bbeta)$ to $m_0(t; \bbeta)$
that $\hat{m}_0(t) \doteq \hat{m}_0(t; \hat{\bbeta}) \rightarrow m_0(t; \bbeta_*) = m_*(t)$ almost surely in $[0, \tau]$. \\
(ii) Write $U(\bbeta_*) \doteq U\{\bbeta_*, \hat{m}_0(t; \bbeta_*)\}$, 
since
\begin{eqnarray*}
& & U\{ \bbeta_*, m_*(\cdot) \}
=  \frac{1}{n} \sum_{i=1}^n  \int_0^\tau \pi_i \bZ_i \big[ m_*(t) dN_i(t) \\[1mm]
& &~~~~~ - Y_i(t) \big\{\exp(-\bbeta_*^{\T} \bZ_i) dt + d m_*(t) \big\} \big],
\end{eqnarray*}
and
\begin{eqnarray*}
& & U\{ \bbeta_*, \hat{m}_0(t; \bbeta_*) \} - U\{ \bbeta_*, m_*(\cdot) \} \\ [2mm]
&=& \frac{1}{n} \sum_{i=1}^n  \int_0^\tau \pi_i \bZ_i \big[ \left\{ \hat{m}_0(t; \bbeta_*)-m_*(t)  \right\} dN_i(t) \\[1mm]
& & ~~~~~ - Y_i(t) d \left\{\hat{m}_0(t; \bbeta_*)-m_*(t) \right\} \big]  \\ [1mm]
&=& \frac{1}{n} \sum_{i=1}^n  \int_0^\tau \pi_i \bZ_i \big[ \left\{ \hat{m}_0(t; \bbeta_*)-m_*(t)  \right\} dN_i(t) \\[1mm]
& & - Y_i(t) \left\{ \frac{ \sum_{i=1}^n \pi_i  dN_i(t)}{ \sum_{i=1}^n \pi_i Y_i(t) } \left\{ \hat{m}_0(t;\bbeta_*) - m_*(t)  \right\}   \right. \\ [1mm]
& &   + \left. \left.  \frac{ \sum_{i=1}^n \pi_i m_*(t) d M_i(t; \bbeta_*, m_*) } { \sum_{i=1}^n \pi_i Y_i(t) }   \right\}   \right] \\ [1mm]
&=&  \frac{1}{n} \sum_{i=1}^n  \int_0^\tau \pi_i [\bZ_i - \bar{\bZ}(t)] d N_i(t) \left\{ \hat{m}_0(t;\bbeta_*) - m_*(t)  \right\}  \\[1mm]
& &   -  \frac{1}{n} \sum_{i=1}^n  \int_0^\tau \pi_i  \bar{\bZ}(t) m_*(t) d M_i(t, \bbeta_*, m_*)\\
&=& \frac{1}{n} \sum_{i=1}^n  \int_0^\tau  \pi_i [\bZ_i - \bar{\bZ}(t)] d N_i(t)
  \bigg\{  - \frac{1}{S_n(t)} \int_t^\tau \frac{S_n(u)}{C_n(u)}   \\ [1mm]
& &~~~~~ \frac{1}{n} \sum_{i=1}^n \pi_i m_*(u) d M_i(u; \bbeta_*, m_*) \bigg\} \\ [1mm]
& &-  \frac{1}{n} \sum_{i=1}^n  \int_0^\tau \pi_i  \bar{\bZ}(t) m_*(t) d M_i(t, \bbeta_*, m_*)\\ [1mm]
&=&  -  \frac{1}{n} \sum_{i=1}^n  \int_0^\tau \pi_i [ \tilde{\bZ}(t) + \bar{\bZ}(t) ] m_*(t) d M_i(t, \bbeta_*, m_*),
\end{eqnarray*}
therefore
\begin{eqnarray*}
& & n^{1/2} U ( \bbeta_* ) \\[1mm]
&=& n^{1/2} U \{ \bbeta_*,\hat{m}_0(t; \bbeta_* ) \} \\ [1mm]
&=& n^{1/2} U\{ \bbeta_*, m_*(\cdot) \} \\ [1mm]
& & + n^{1/2} [ U\{ \bbeta_*, \hat{m}_0(t; \bbeta_*) \} - U\{ \bbeta_*, m_*(\cdot) \} ] \\ [1mm]
&=&  n^{-1/2}  \sum_{i=1}^n  \int_0^\tau \pi_i \left[ \bZ_i - \bar{\bZ}(t)- \tilde{\bZ}(t) \right] \\ [1mm]
& & ~~~~~ m_*(t) d M_i(t, \bbeta_*, m_*) \\ [1mm]
&=& n^{-1/2}  \sum_{i=1}^n  \int_0^\tau \pi_i \left[ \bZ_i - \mu_{\bz}(t)- \tilde{\mu}_{\bz}(t) \right] \\[1mm]
& & ~~~~~ m_*(t) d M_i(t, \bbeta_*, m_*)+o_p(1) \\[1mm]
&=&  n^{-1/2}  \sum_{i=1}^n  \int_0^\tau \left[ \bZ_i - \mu_{\bz}(t)- \tilde{\mu}_{\bz}(t) \right] m_*(t) d M_i(t, \bbeta_*, m_*)\\[1mm]
& & + n^{-1/2}  \sum_{i=1}^n  \int_0^\tau (\pi_i-1) \left[ \bZ_i - \mu_{\bz}(t)- \tilde{\mu}_{\bz}(t) \right] \\[1mm]
& & ~~~~~ m_*(t) d M_i(t, \bbeta_*, m_*)+o_p(1) \\[1mm]
&=&  n^{-1/2}  \sum_{i=1}^n  \int_0^\tau \left[ \bZ_i - \mu_{\bz}(t)- \tilde{\mu}_{\bz}(t) \right] m_*(t) d M_i(t, \bbeta_*, m_*) \\[1mm]
& & - n^{-1/2}  \sum_{i=1}^n  \int_0^\tau (1-\delta_i)(1-\xi_i/p) \left[ \bZ_i - \mu_{\bz}(t)- \tilde{\mu}_{\bz}(t) \right] \\[1mm]
& & ~~~~~ m_*(t) d M_i(t, \bbeta_*, m_*)+o_p(1).
\end{eqnarray*}
Let $\mathcal{F}_i$ be the $\sigma$-field generated by $\{\tilde{T}_i, \delta_i, \bZ_i\}$, and
$$
\eta_i = (1-\delta_i)  \int_0^\tau \left[ \bZ_i - \mu_{\bz}(t)- \tilde{\mu}_{\bz}(t) \right] m_*(t) d M_i(t, \bbeta_*, m_*).
$$
It is easy to see that $E(1-\xi_i/p|\mathcal{F}_i) = 0$, $E\{\eta_i (1-\xi_i/p) | \mathcal{F}_i  \} = E\{ \eta_i E( 1-\xi_i/p |\mathcal{F}_i) \}=0$.
Following the arguments in \cite{Kulich2000} and \cite{Lu2006}, we have $ \mbox{var} \{ \eta_i (1-\xi_i/p) \} = (1-p)/p \cdot \{ E[\eta_i^{\otimes 2}] - E[\eta_i]^{\otimes 2} \}  = \Sigma_2$.
Also, conditional on $\mathcal{F}_i$, $\{ \eta_i(1-\xi_i/p), i=1, \ldots, n \}$ and the first term of $n^{1/2} U ( \bbeta_* )$ are uncorrelated, and hence $n^{1/2} U(\boldsymbol{\beta}_*)$ is asymptotically normal with mean zero and variance-covariance matrix $\Sigma=\Sigma_1+\Sigma_2.$
Thus, the Taylor expansion of $U(\hat{\bbeta})$ at $\bbeta_*$ gives
\begin{eqnarray*}
& & n^{1/2} (\hat{\bbeta}-\bbeta_*)   \\ [1mm]
&=& - A^{-1} n^{1/2} U(\bbeta_*) + o_p(1)   \\ [1mm]
&=&   - A^{-1} n^{-1/2} \sum_{i=1}^n \int_0^\tau \pi_i \left[ \bZ_i - \mu_{\bz}(t)- \tilde{\mu}_{\bz}(t) \right] \\[1mm]
& &~~~~~ m_*(t) d M_i(t, \bbeta_*, m_*)+o_p(1).
\end{eqnarray*}
(iii) Let $\zeta_i = \int_0^\tau \left[ \bZ_i - \mu_{\bz}(t)- \tilde{\mu}_{\bz}(t) \right] m_*(t) d M_i(t, \bbeta_*, m_*)$.
Write
\begin{eqnarray*}
& & n^{1/2} \{ \hat{m}_0(t) - m_*(t) \}   \\ [1mm]
&=& n^{1/2} \{ \hat{m}_0(t; \hat{\bbeta}) - \hat{m}_0(t; \bbeta_*) \}
         + n^{1/2} \{ \hat{m}_0(t; \bbeta_*) - m_*(t) \}  \\ [1mm]
&=& \left( \left. \frac{\partial \hat{m}_0(t;\boldsymbol{\beta})}{\partial \boldsymbol{\beta}} \right|_{\boldsymbol{\beta}=\boldsymbol{\beta}_*}  \right)^{\T}
     n^{1/2} (\hat{\bbeta} - \bbeta_*)  \\[1mm]
& & + n^{1/2} \{ \hat{m}_0(t; \bbeta_*) - m_*(t) \} + o_p(1).
\end{eqnarray*}

Now, we establish the formula of $\hat{m}_0(t; \bbeta_*)-m_*(t)$.
By inserting $\{ \bbeta_*, m_*(t)\}$ and
$\{ \bbeta_*, \hat{m}_0(t; \bbeta_*) \}$ into $(\ref{eq113})$, respectively, we have
\begin{eqnarray}
 & & \sum_{i=1}^n \pi_i \big[ m_*(t) dN_i(t)-
    Y_i(t) \big\{\exp(-\bbeta_*^{\T} \bZ_i) dt    \nonumber \\[1mm]
 & & ~~~~~ + d m_*(t) \big\} \big] = \sum_{i=1}^n \pi_i m_*(t) d M_i(t; \bbeta_*, m_*),  \label{eqap1}  \quad \quad \\  [1mm]
 & & \sum_{i=1}^n \pi_i \big[ \hat{m}_0(t;\bbeta_*) dN_i(t)   \nonumber   \\ [1mm]
 & & ~~~~~ - Y_i(t) \big\{\exp(-\bbeta_*^{\T} \bZ_i) dt + d \hat{m}_0(t;\bbeta_*) \big\} \big] = 0.
\label{eqap2}
\end{eqnarray}
Subtracting $(\ref{eqap1})$ from $(\ref{eqap2})$,
\begin{eqnarray*}
& & \sum_{i=1}^n \pi_i \big[ \big\{ \hat{m}_0(t;\bbeta_*) - m_*(t)  \big\} dN_i(t) \\ [1mm]
& & ~~~~~ - Y_i(t)  d  \big\{ \hat{m}_0(t;\bbeta_*) - m_*(t)  \big\} \big]   \\[1mm]
&=&  - \sum_{i=1}^n \pi_i m_*(t) d M_i(t; \bbeta_*, m_*),
\end{eqnarray*}
which is equivalent to
\begin{eqnarray*}
& & \frac{ \sum_{i=1}^n \pi_i  dN_i(t)}{ \sum_{i=1}^n \pi_i Y_i(t) }  \left\{ \hat{m}_0(t;\bbeta_*) - m_*(t)  \right\}  \\ [1mm]
& & ~~~~~ - d  \left\{ \hat{m}_0(t;\bbeta_*) - m_*(t)  \right\}   \\[1mm]
&=&  - \frac{ \sum_{i=1}^n \pi_i m_*(t) d M_i(t; \bbeta_*, m_*) } { \sum_{i=1}^n \pi_i Y_i(t) }.
\end{eqnarray*}
Then
\begin{eqnarray*}
& & \hat{m}_0(t;\bbeta_*) - m_*(t)    \\[1mm]
&=& - \frac{1}{S_n(t)} \int_t^\tau \frac{S_n(u)}{C_n(u)}
 \frac{1}{n} \sum_{i=1}^n \pi_i m_*(t) d M_i(t; \bbeta_*, m_*).
\end{eqnarray*}

Hence
\begin{eqnarray*}
& & n^{1/2} \{ \hat{m}_0(t) - m_*(t) \}   \\ [1mm]
&=& m_*(t) \mu_\bz(t)^{\T} A^{-1} n^{-1/2} \sum_{i=1}^n \pi_i \zeta_i
        - \frac{1}{S_n(t)} \int_t^\tau \frac{S_n(u)}{C_n(u)}    \\ [1mm]
& & ~~~~~ n^{-1/2} \sum_{i=1}^n \pi_i m_*(u) d M_i(u; \bbeta_*, m_*) + o_p(1)  \\ [1mm]
&=& m_*(t) \mu_\bz(t)^{\T} A^{-1} n^{-1/2} \sum_{i=1}^n \pi_i \zeta_i
        - \frac{1}{S(t)} \int_t^\tau \frac{S(u)}{C(u)} \\ [1mm]
& & ~~~~~ n^{-1/2} \sum_{i=1}^n \pi_i m_*(u) d M_i(u; \bbeta_*, m_*) + o_p(1)  \\ [1mm]
&=& n^{-1/2} \sum_{i=1}^n \pi_i \bigg[ m_*(t) \mu_\bz(t)^{\T} A^{-1} \zeta_i
                  - \frac{1}{S(t)} \int_t^\tau \frac{S(u)}{C(u)}  \\ [1mm]
& & ~~~~~ m_*(u) d M_i(u; \bbeta_*, m_*)   \bigg] + o_p(1)  \\ [1mm]
&=& n^{-1/2} \sum_{i=1}^n \varphi_i(t) - n^{-1/2} \sum_{i=1}^n (1-\delta_i)(1-\xi_i/p) \varphi_i(t)  \\[1mm]
& & ~~~~~  + o_p(1),
\end{eqnarray*}
where
\begin{eqnarray*}
\varphi_i(t) &=& m_*(t) \mu_\bz(t)^{\T} A^{-1} \zeta_i  \\ [1mm]
  & & - \frac{1}{S(t)} \int_t^\tau \frac{S(u)}{C(u)}  m_*(u) d M_i(u; \bbeta_*, m_*),
\end{eqnarray*}
$S(t)$ and $C(t)$ are the uniform limit of $S_n(t)$ and $C_n(t)$ respectively.
Because $\{\varphi_i(t), i=1, \ldots, n)\}$ are independent mean zero random variables for each $t$,
by \cite{pollard1990empirical},
 $n^{1/2} \{ \hat{m}_0(t) - m_*(t) \} (0 \leq t \leq \tau)$ converges weakly to a mean zero Gaussian process,
 whose covariance function at $(s, t)$ is
\begin{eqnarray*}
 & & E[\varphi_1(s) \varphi_1(t)] + \frac{1-p}{p} E[(1-\delta)\varphi_1(s) \varphi_1(t)]  \\ [1mm]
 & & ~~~ - \frac{1-p}{p} E[(1-\delta)\varphi_1(s)] E[(1-\delta)\varphi_1(t)].
\end{eqnarray*}
This ends the proof.

\section*{Acknowledgements}

The authors wish to thank the Editor Yuedong Wang, the Associate Editor and two anonymous referees for their valuable and helpful comments.
The authors also wish to express their appreciation to Professor Ying Qing Chen and
Cindy Zhang for their invaluable assistance of the original version.
Ma's work is partially supported by National Institutes of Health grant R01 HL113548.
Shi's work is supported by Natural Science Foundation of Fujian Province, China (2016J01026), the Institute of Meteorological Big Data-Digital Fujian and Fujian Key Laboratory of Data Science and Statistics, China.  Zhou's work is supported by the State Key Program of National Natural Science Foundation of China (71331006), the State Key Program in the Major Research Plan of National Natural Science Foundation of China (91546202).

\end{document}